\documentclass[12pt]{article}
\usepackage{amssymb}
\usepackage{latexsym}
\usepackage{amsmath}

\newtheorem{theo}{Theorem}[section]

\newtheorem{lemma}[theo]{Lemma}
\newtheorem{claim}[theo]{Claim}

\newtheorem{con}[theo]{Conjecture}

\newtheorem{fact}[theo]{Fact}

\def\q{\hspace*{\fill}$\Box$\medskip}
\begin{document}

\title{Bipartite Ramsey numbers of cycles}
\author{Zilong Yan \thanks{School of Mathematics, Hunan University, Changsha 410082, P.R. China. Email: zilongyan@hnu.edu.cn.} \and Yuejian Peng \thanks{ Corresponding author. School of Mathematics, Hunan University, Changsha, 410082, P.R. China. Email: ypeng1@hnu.edu.cn. \ Supported in part by National Natural Science Foundation of China (No. 11931002).}
}

\maketitle

\begin{abstract}
Given bipartite graphs $H_1$, \dots , $H_k$, the bipartite Ramsey number $br(H_1,\dots, H_k)$ is the minimum integer $N$ such that any $k$-edge-coloring of complete bipartite graph $K_{N, N}$ contains a monochromatic $H_i$ in color $i$ for some $i\in \{1, 2, \ldots, k\}$. There are considerable results on asymptotic values of bipartite Ramsey numbers of cycles. For exact values, Zhang-Sun \cite{Zhangs} determined $br(C_4, C_{2n})$, Zhang-Sun-Wu \cite{Zhangsw} determined $br(C_6, C_{2n})$, and Gholami-Rowshan \cite{GR} determined $br(C_8, C_{2n})$. In this paper, we determine completely the exact values of $br(C_{2n}, C_{2m})$ for all $n\ge m\ge 5$, this answers a question concerned by Buci\'c-Letzter-Sudakov \cite{BLS}, Gholami-Rowshan \cite{GR}, and Zhang-Sun \cite{Zhangs}. Some observations obtained in the proof have their own interests (in our opinion).
\medskip

\end{abstract}

\section{Introduction}
Let $H_1,\dots , H_k$ be graphs. The Ramsey number $r(H_1, \dots, H_k)$ is the minimum integer $N$ such that any $k$-edge-coloring of $K_N$ contains a monochromatic $H_i$ in color $i$ for some $i\in \{1, 2, \ldots, k\}$. If $H_1=H_2=\dots=H_k=H$, then simplify it as $r^k(H)$. Let $H_1$, \dots , $H_k$ be bipartite graphs. The bipartite Ramsey number $br(H_1,\dots, H_k)$ is the minimum integer $N$ such that any $k$-edge-coloring of complete bipartite graph $K_{N, N}$ contains a monochromatic $H_i$ in color $i$ for some $i\in \{1, 2, \ldots, k\}$. If $H_1=H_2=\dots=H_k=H$, then simplify it as $br^k(H)$. The existence of such an integer $N$ is guaranteed by a result of Ramsey \cite{Ramsey} and Erd\H os and Rado \cite{ES}. The study of bipartite Ramsey number was initiated in the early 1970s by Faudree and Schelp \cite{FS}, and Gy\'arf\'as and Lehel \cite{GL}. Let $P_n$ denote a path with $n$ vertices. They determined the bipartite ramsey numbers of paths. 
Bipartite Ramsey numbers are also studied for the complete bipartite graphs, and the first to consider this is Beineke and Schwenk \cite{BS} in 1976. Similar to the case of ordinary Ramsey numbers, the best known lower bound on $br^k(K_{n, n})$, due to Hattingh and Henning \cite{HH} and the best known upper bound, due to Conlon \cite{Con}, are still exponentially apart.

There are considerable nice results on Ramsey numbers of cycles (see \cite{BE, FS1, Rosta, Luczak, FL, BS2, JS}). In this paper, we focus on bipartite Ramsey numbers of cycles. Let $C_n$ denote a cycle with $n$ vertices.  Goddard, Henning and Oellermann \cite{GHO} showed that $br^3(C_4)=11$. Joubert \cite{Jou} showed that
$$br(C_{2t_1}, C_{2t_2}, \dots, C_{2t_k})\le k(t_1+t_2+\dots+t_k-k+1),$$
where $t_i$ is an integer and $2\le t_i\le 4$ for all $1\le i\le k$. 
Desiasio, Gy\'arf\'as, Krueger, Ruszink\'o and S\'ark\H ozy \cite{DGKRS} showed that $br^k(C_{2n})\ge (2k-4)n$ for $k\ge 5$ and $br^4(C_{2n})\ge 5n$. Buci\'c, Letzter and Sudakov showed that $br^3(C_{2n})=(3+o(1))n$ in \cite{BLS}, and $br^k(C_{2n})\le (2k-3+o(1))n$ for $k\ge 5$ and $br^4(C_{2n})=(5+o(1))n$ in \cite{BLS1}. Liu and Peng \cite{LiuP} gave the asymptotic value of  $br(C_{2\lfloor\alpha_1n\rfloor}, \dots, C_{2\lfloor\alpha_rn\rfloor})$ when $r\ge 3$, $\alpha_1, \alpha_2>0$ and $\alpha_{j+2}\ge [(j+2)!-1]\sum_{i=1}^{j+1}\alpha_i$ for $1\le j\le r-2$. Luo and Peng \cite{LP} gave the asymptotic value of $br(C_{2\lfloor\alpha_1n\rfloor}, C_{2\lfloor\alpha_2n\rfloor}, C_{2\lfloor\alpha_3n\rfloor})$ for any $\alpha_1, \alpha_2, \alpha_3>0$. Several of the above mentioned asymptotic results used a method initiated  by {\L}uczak which applies Szemer\'edi's Regularity Lemma  to  obtain the asymptotic values of (bipartite) ramsey numbers of cycles  by showing  the existence of monochromatic connected matchings in  almost (bipartite) complete graphs (the reduced graph guaranteed by Regularity Lemma). Indeed, Letzter\cite{letz} showed that these can be
further reduced to problems about ramsey numbers of  monochromatic connected matchings. What about the exact values of bipartite ramsey numbers of cycles themselves?
Combining the results by Beineke and Schwenk \cite{BS}, Zhang and Sun \cite{Zhangs},  Zhang, Sun and Wu \cite{Zhangsw}, and Gholami and Rowshan \cite{GR}, we have
\begin{theo}(\cite{BS}, \cite{GR}, \cite{Zhangs},\cite{Zhangsw})\label{le4}
$$ br(C_{2n}, C_{2m})=\left\{
\begin{array}{rcl}
5   &  &m=2; \ n=2 \ or \ 3, \\
n+1   &  &m=2; \ n\ge 4,\\
6   &  &m=3; \ n=3, \\
n+2   &  &m=3; \ n\ge 4,\\
8   &  &m=4; \ n=4, \\
n+3   &  &m=4; \ n\ge 5.
\end{array}\right. $$
\end{theo}
Zhang, Sun and Wu \cite{Zhangsw}, and Gholami and Rowshan \cite{GR} conjectured that
\begin{con}(\cite{Zhangsw}, \cite{GR})\label{conj}
$br(C_{2n}, C_{2m})=m+n-1$ for $n>m$, and $br(C_{2n}, C_{2m})=m+n$ for $n=m$.
\end{con}

{\em The Lower bound of Conjecture \ref{conj}} (due to Zhang-Sun-Wu \cite{Zhangsw}):

If $n>m$, let $G=G(X_1\cup X_2, Y)$ be a complete bipartite graph with bipartition $(X_1\cup X_2)\cup Y$, where $|X_1|=m-1$, $|X_2|=n-1$, and $|Y|=n+m-2$. Color the edges in $G(X_1, Y)$ red and the edges in $G(X_2, Y)$ blue. Then there is neither a red $C_{2m}$ nor a blue $C_{2n}$ in $G$. So $br(C_{2n}, C_{2m})\ge m+n-1$.

If $m=n$, let $G=G(X_1\cup X_2\cup\{x\}, Y_1\cup Y_2\cup\{y\})$ be a complete bipartite graph with bipartition $(X_1\cup X_2\cup\{x\})\cup(Y_1\cup Y_2\cup\{y\})$, where $|X_1|=|X_2|=|Y_1|=|Y_2|=m-1$. Color the edges in $G(X_1, Y_1)\cup G(X_2, Y_2)\cup G(\{x\}, Y_1\cup Y_2\cup\{y\})$ red and the edges in $G(X_1, Y_2)\cup G(X_2, Y_1)\cup G(X_1\cup X_2, \{y\})$ blue. Then there is neither a red $C_{2m}$ nor a blue $C_{2m}$ in $G$. So $br(C_{2n}, C_{2m})\ge m+n$.
\q

In this paper, we confirm Conjecture \ref{conj} for all the remaining cases and obtain the following results.
\begin{theo}\label{main}
$$ br(C_{2n}, C_{2m})=\left\{
\begin{array}{rcl}
n+m-1   &  &n\neq m, \ n, m\ge 5, \\
m+n   &  &n=m, \ n, m\ge 5.
\end{array}\right. $$
\end{theo}

In Section \ref{sec2}, we will give the proof of Theorem \ref{main},  but leave the proofs of several key lemmas in Section \ref{sec3}.

{\bf Notations.} For a graph $G$, let $V(G)$ and $E(G)$ denote its vertex set and edge set respectively. Given a vertex $x$ in $G$, $d_G(x)$, for short $d(x)$, denotes its degree in $G$.
Given a subset $U$ of $V(G)$, $G[U]$ denotes the subgraph of $G$ induced by $U$, and $E(G[U])$, for short $E(U)$, denotes the edge set of $G[U]$. Let $e(U)=|E(U)|$. Given a subset $U$ of $V(G)$ and a vertex $x$ in $G$, $d_U(x)$ denotes the number of edges with one endpoint $x$ and one endpoint in $U$. Given a red/blue-edge-coloring of $G$, $U\subseteq V(G)$ and $x\in V(G)$, $d_U^r(x)$ ($d_U^b(x)$) denotes the number of red (blue) edges with one endpoint $x$ and one endpoint in $U$.
Let $G(X, Y)$, sometimes for short $(X, Y)$, be a bipartite graph with bipartition $X\cup Y$. Given any red/blue-edge-coloring of $G$, the red (blue) edge set is denoted by $E^r(G)$ ($E^b(G)$), and the number of red (blue) edges is denoted by $e^r(G)$ ($e^b(G)$). Given $A\subseteq V(G)$ and $B\subseteq V(G)$, $E(A, B)$ denotes all edges with one endpoint in $A$ and one endpoint in $B$. Let $e(A, B)=|E(A, B)|$. And $E^r(A, B)$ ($E^b(A, B)$) denotes all red (blue) edges with one endpoint in $A$ and one endpoint in $B$. Let $e^r(A, B)=|E^r(A, B)|$ and $e^b(A, B)=|E^b(A, B)|$.
Let $C$ be a cycle. Assuming that we walk along the cycle clockwisely, we write $v^-$ ($v^+$) for the predecessor (successor) of $v\in V(C)$. Given $U\subseteq V(C)$, $U^+=\{u^+| u\in U\}$ and $U^-=\{u^-| u\in U\}$. The number of vertices in a cycle $C$ is denoted by $|C|$, and the number of vertices in a path $P$ is denoted by $|P|$.

\section{Proof of Theorem \ref{main}}\label{sec2}

The Lower bound is  given after Conjecture \ref{conj}, we only need to show the upper bound. Let us give a very rough sketch of the proof. We apply induction on $n+m$.  Suppose on the contrary that there is a red/blue-edge-coloring of $G(X, Y)=K_{N, N}$ such that there is neither a blue $C_{2n}$ nor a red $C_{2m}$. By induction hypothesis,  $G$ has a blue cycle $C$ with $|C|=2(n-1)$. Indeed we will show that the length of a longest blue cycle is $2(n-1)$ (Lemma \ref{c2}, a crucial lemma). Let $A=X\setminus V(C)$ and $B=Y\setminus V(C)$. We will show that  there exists a `large' blue component $(A_1, B_1)$ in $G-V(C)$ with $A_1\subseteq A$ and $B_1\subseteq B$ and $|A_1|+|B_1|\ge |A|$ (See Lemma \ref{c5}). We extend the blue component $(A_1, B_1)$ in $G-V(C)$ to the blue cycle $C$, namely, define
$$X_1=\{u\in V(C)|uv\in E^b(G) \ for \ some \ v\in B_1\},$$
$$Y_1=\{u\in V(C)|uv\in E^b(G) \ for \ some \ v\in A_1\}.$$
Let $X_2=X\setminus(A\cup X_1)$, $Y_2=Y\setminus(B\cup Y_1)$, $A_2=A\setminus A_1$ and $B_2=B\setminus B_1$. Therefore we get two red complete bipartite graphs $(A_1, B_2\cup Y_2)$ and $(B_1, A_2\cup X_2)$. We will show that there is a red $C_{2m}$ by taking advantage of these two red complete bipartite graphs. There are red paths in these two red complete bipartite graphs, but how to glue two red paths to obtain a red $C_{2m}$? It is easy to show that if there are two disjoint red edges between these two red complete bipartite graphs and these two red complete bipartite graphs have enough number of vertices, then we are fine (Lemma \ref{c3}). How to guarantee that there are two disjoint red edges between these two red complete bipartite graphs and the number of vertices in these two red complete bipartite graphs are enough? These are the main technical part we need to deal with. If both $X_1$ and $Y_1$ are non-empty, then using the fact that $C_{2(n-1)}$ is a longest blue cycle (guaranteed by Lemma \ref{c2}, to be proved in Section \ref{sec3}) and its implications (Fact \ref{c7}), we have two disjoint red edges in $(X_2, Y_2)$; by some quantitative analysis, we can also show numbers of vertices in two red complete bipartite graphs are enough, this is a relatively easy case. For another case,   we will show a crucial lemma (Lemma \ref{c1}) to  guarantee  that there is no blue $K_{n-1, n-1}$ in $G$. However, if there are not two disjoint red edges, we may find a blue $K_{n-1, n-1}$ through some analysis (Lemma \ref{c4} to be proved in Section \ref{sec3} is applied sometimes), and get a contradiction. This is a very rough sketch, to make everything fit, we need more subtle arguments. The precise argument will be given in this section. Before giving the precise proof, we state the above mentioned lemmas  and leave  proofs of some of them  in Section\ref{sec3}. In the  proofs of some lemmas, we obtain some observations having their own interests (in our own opinion) in Section \ref{sec3}.

\begin{lemma}\label{c3}
Let $(A_1, B_1)$ and $(A_2, B_2)$ be two complete bipartite graphs, where $A_1, B_1, A_2, B_2$ are pairwisely disjoint. If there is one edge in $(A_1, A_2)$ and one edge in $(B_1, B_2)$, then there are all cycles with length from 6 to $2(min\{|A_1|, |B_1|\}+min\{|A_2|, |B_2|\})$. If there are two disjoint edges in $(A_1, A_2)$, then there are all cycles with length from 6 to $min\{2|A_1|-1, 2|B_1|+1\}+min\{2|A_2|-1, 2|B_2|+1\}$.
\end{lemma}
{\em Proof of Lemma \ref{c3}.} If there is one edge in $(A_1, A_2)$ and one edge in $(B_1, B_2)$. Without loss of generality, assume that $a_1a_2$ and $b_1b_2$ are edges, where $a_1\in A_1$, $a_2\in A_2$, $b_1\in B_1$ and $b_2\in B_2$. Since $(A_1, B_1)$ is a complete bipartite graph, we can find all paths with endpoints $a_1$ and $b_1$ of order 2 to $2min\{|A_1|, |B_1|\}$ in $(A_1, B_1)$. Similarly, we can find all paths with endpoints $a_2$ and $b_2$ of order 2 to $2min\{|A_2|, |B_2|\}$. These paths and $a_1a_2$ and $b_2b_2$ form all cycles with length from 6 to $2(min\{|A_1|, |B_1|\}+min\{|A_2|, |B_2|\})$. Similarly, if there are two disjoint edges in $(A_1, A_2)$, then there are all cycles with length from 6 to $min\{2|A_1|-1, 2|B_1|+1\}+min\{2|A_2|-1, 2|B_2|+1\}$.
\q

\begin{lemma}\label{c5}
Let $G(A, B)$ be a bipartite graph with $|A|=|B|\ge m$. If a red/blue-edge-coloring of $G$ has no red $C_{2m}$, then there is a blue component $(A_1, B_1)$ with $A_1\subset A$ and $B_1\subset B$ such that $|A_1|+|B_1|\ge |A|$.
\end{lemma}

In the proof of Lemma \ref{c5}, we will apply Lemma \ref{c3} and the following results by Moon-Moser \cite{MM}.
\begin{theo}\label{n1}(\cite{MM})
Let $G(X, Y)$ be a bipartite graph with $|X|=|Y|=n$ and $2n\ge 8$. If $d(x)+d(y)\ge n+1$ for every pair of nonadjacent vertices $x\in X$ and $y\in Y$, then $G$ has a hamiltonian cycle.
\end{theo}

{\em Proof of Lemma \ref{c5}.} Assume that $(A_1, B_1)$ is a largest blue component and $|A_1|+|B_1|<|A|$. Choose an $m$-set in $A$ and an $m$-set in $B$. Since there is no red $C_{2m}$ in the two $m$-sets, by Theorem \ref{n1}, there is a blue edge $ab$ such that the red neighbors of $a$ and $b$ is less than $m+1$. So there is a blue component $(A_1', B_1')$ and $|A_1'|+|B_1'|\ge m$. By the maximality of  $(A_1, B_1)$, we have $|A_1|+|B_1|\ge m$. Let $A_2= A\setminus A_1$ and $B_2= B\setminus B_1$. So $G(A_1, B_2)$ and $G(A_2, B_1)$ are red complete bipartite graphs with $|A_2|>|B_1|$ and $|B_2|> |A_1|$. By Lemma \ref{c3}, there is no disjoint red edges in $G(A_2, B_2)$, otherwise there is a red $C_{2m}$. So there is a blue $K_{|A_2|-1, |B_2|}$ or a blue $K_{|A_2|, |B_2|-1}$ in $G(A_2, B_2)$. Such a blue component has at least $|A_2|+|B_2|-1\ge 2|A|-|A_1|-|B_1|-1\ge |A|$ vertices, a contradiction.\q

\begin{lemma}\label{c1}
Let $N=m+n$ for $m=n\ge 5$ and $N=m+n-1$ for $n\ge m+1\ge 5$. If there is a red/blue-edge-coloring of $G(X, Y)=K_{N, N}$ such that there is neither a blue $C_{2n}$ nor a red $C_{2m}$. Then there is no blue $K_{n-1, n-1}$ in $G$.
\end{lemma}

\begin{lemma}\label{c2}
Let $N=m+n$ for $m=n$ and $N=m+n-1$ for $n\ge m+1$. Assume that there is a red/blue-edge-coloring of $G(X, Y)=K_{N, N}$ such that there is neither a blue $C_{2n}$ nor a red $C_{2m}$. Assume that $C$ is a longest blue cycle with $|V(C)|=2c$ and $C'$ is a longest red cycle with $|V(C')|=2t$. If $n\ge m+1$, then $c\le n-1$. If $n=m$, then $c\le m-1$ and $t\le m-1$.
\end{lemma}

Proofs of Lemmas \ref{c1} and \ref{c2} will be given in Section \ref{sec3}. The following easy fact will be used in the proof of Theorem \ref{main}.

\begin{fact}\label{c7}
Let $G$ be a graph. Let $C$ be a longest cycle in $G$. Assume that $P\subset V(G)\setminus V(C)$ is a path with endpoints $x$ and $y$. Let $ux, vy\in E(G)$ and $u, v\in V(C)$. Then $u^-v^-\notin E(G)$, $u^+v^+\notin E(G)$, $u^+v^{++}\notin E(G)$, $u^{--}v^{-}\notin E(G)$, $u^{--}y\notin E(G)$, and $v^{++}x\notin E(G)$.
\end{fact}
{\em Proof of Fact \ref{c7}.} If $u^-v^-\in E(G)$, then  $v$ (along $C$ clockwisely)$\rightarrow u^- \rightarrow v^-$ (along $C$ counterclockwisely)$\rightarrow u \rightarrow x\stackrel{P}{\longrightarrow}y\rightarrow v$ forms a cycle with length $|C|+|P|>|C|$, a contradiction to that $C$ is a longest cycle in $G$. Similarly, $u^+v^+\notin E(G)$, $u^+v^{++}\notin E(G)$ and $u^{--}v^{-}\notin E(G)$. If $u^{--}y\in E(G)$, then  $u$ (along $C$ clockwisely)$\rightarrow u^{--}\rightarrow y$ $\stackrel{P}{\longrightarrow}x\rightarrow u$ forms a cycle with length $|C|+|P|-1>|C|$, a contradiction to that $C$ is a longest cycle in $G$. Similarly, $v^{++}x\notin E(G)$.\q

{\em Proof of Theorem \ref{main}.} By the Lower bound given after Conjecture \ref{conj}, we only need to show the upper bound. Without loss of generality, we assume that $n\ge m$.  Let $N=m+n$ for $m=n$ and $N=m+n-1$ for $n\ge m+1$. Using induction on $n+m$. By Theorem \ref{le4}, the conclusion holds for $br(C_8, C_{2n})$, where $n\ge 4$. We assume that the conclusion holds for $br(C_{2k}, C_{2l})$ provided $k, l\ge 4$ and $k+l\le m+n-1$. We aim to show the conclusion for $br(C_{2n}, C_{2m})$ if $n\ge m\ge 5$. Suppose that there is a red/blue-edge-coloring of $G(X, Y)=K_{N, N}$ such that there is neither a blue $C_{2n}$ nor a red $C_{2m}$.  In view of Lemma \ref{c2}, we may assume that the length of a longest blue cycle is at most $2(n-1)$. By induction hypothesis, $br(C_{2(n-1)}, C_{2m})\le N$, so $G$ has a blue cycle $C$ with $|C|=2(n-1)$. Let $A=X\setminus V(C)$ and $B=Y\setminus V(C)$. Since $|A|=|B|\ge m$ and $G$ contains no red $C_{2m}$, by Lemma \ref{c5}, there exists a blue component $(A_1, B_1)$ with $A_1\subseteq A$ and $B_1\subseteq B$ and $|A_1|+|B_1|\ge |A|$. Let
$$X_1=\{u\in V(C)|uv\in E^b(G) \ for \ some \ v\in B_1\},$$
$$Y_1=\{u\in V(C)|uv\in E^b(G) \ for \ some \ v\in A_1\}.$$
Let $X_2=X\setminus(A\cup X_1)$, $Y_2=Y\setminus(B\cup Y_1)$, $A_2=A\setminus A_1$ and $B_2=B\setminus B_1$. Therefore $(A_1, B_2\cup Y_2)$ and $(B_1, A_2\cup X_2)$ are red complete bipartite graphs. Note that any vertex in $A_1$ is connected to any vertex in $B_1$ by a blue path. Therefore for any $u\in Y_1$, $u^-, u^+\notin X_1$, (recall that $u^-, u^+$ are predecessor and successor of $u$ along $C$ clockwisely) otherwise there is a longer blue cycle, a contradiction to that $C$ is a longest blue cycle.

{\em Case 1.} $|X_1|, |Y_1|\ge 1$.

Recall that $Y_1^+=\{y^+\in V(C)| y\in Y_1\}$ and $Y_1^-=\{y^-\in V(C)| y\in Y_1\}$. Since $Y_1^+\cup Y_1^-\subseteq X_2$, we have $|X_1|+|Y_1^+\cup Y_1^-|\le n-1$, and $|X_1|+|Y_1|+1\le |X_1|+|Y_1^+\cup Y_1^-|\le n-1$, i.e. $|X_1|+|Y_1|\le n-2$. So $|X_2|+|Y_2|\ge n$. By Fact \ref{c7}, there exist two disjoint red edges $u^+v^+, u^-v^-\in (X_2, Y_2)$ for $u\in X_1$ and $v\in Y_1$. Recall that $|A_1|+|B_1|\ge |A|$, then $|A_1|\ge |B_2|$ and $|B_1|\ge |A_2|$. Note that $(A_1, B_2\cup Y_2)$ and $(B_1, A_2\cup X_2)$ are red complete bipartite graphs. If $|X_2|> |B_1|-|A_2|$, then by Lemma \ref{c3}, $G$ contains all red cycles of length from $6$ to $2(|B_1|+|B_2|)$, a contradiction to that $G$ has no red $C_{2m}$. Similarly, if $|Y_2|> |A_1|-|B_2|$, then there exists a red cycle of size $2(|A_1|+|A_2|)$, a contradiction to that $G$ has no red $C_{2m}$.

If $|X_2|\le |B_1|-|A_2|$ and $|Y_2|\le |A_1|-|B_2|$, then by Lemma \ref{c3}, $G$ contains all red cycles of length  from $6$ to $2(|X_2|+|Y_2|+|A_2|+|B_2|)-2$. Recall that $|X_2|+|Y_2|\ge n\ge m$ and $G$ contains no red $C_{2m}$. Then $2(|X_2|+|Y_2|)=2n=2m$ and $|A_2|=|B_2|=0$, i.e. $|A_1|=|B_1|=|A|$. If there is a red edge in $(A_1, B_1)$, then by Lemma \ref{c3}, $G$ contains all red cycles of length from $6$ to  $2(|X_2|+|Y_2|)$, a contradiction to that $G$ contains no red $C_{2m}$. So $(A_1, B_1)$ forms a blue complete bipartite graph, then $(A_1, B_1)$ contains a blue cycle of length $2n$, a contradiction.

{\em Case 2.} $Y_1=\emptyset$ and $X_1=\emptyset$.

By Lemma \ref{c1}, $G$ has no blue $K_{n-1, n-1}$, then there is a red edge in $(X_2, Y_2)$. Recall that $(A_1, B_2\cup Y_2)$ and $(B_1, A_2\cup X_2)$ are red complete bipartite graphs. Since $G$ has no red $C_{2m}$, by Lemma \ref{c3}, there is no red edges in $(A_1, B_1)$, i.e. $(A_1, B_1)$ is a blue complete bipartite graph. If $|A_1|\le m-1$, then $|A_2|\ge 1$, and consequently $|B_1|\le m-1$ since otherwise the red complete bipartite graph $(B_1, X_2\cup A_2)$ contains a red $C_{2m}$. Therefore $|X_2\cup A_2|, |Y_2\cup B_2|\ge n.$ ($|X_2\cup A_2|, |Y_2\cup B_2|\ge m+1$ if $n=m$.) By Lemma \ref{c3}, there is no disjoint red edges in $(X_2\cup A_2, Y_2\cup B_2)$. Then there is a blue $K_{n-1, n-1}$ in $(X_2\cup A_2, Y_2\cup B_2)$, a contradiction to Lemma \ref{c1}. If $|A_1|=|B_1|\ge m$, recall that we proved that $(A_1, B_1)$ is blue complete bipartite graph and $G$ contains no blue $C_{2n}$, then $n\ge m+1$. Therefore $(A_1, Y_2)$ contains a red $C_{2m}$, a contradiction.

{\em Case 3.} $X_1\neq\emptyset$ and $Y_1=\emptyset$. (By symmetry, the case $X_1=\emptyset$ and $Y_1\neq\emptyset$ is similar.)

\begin{claim}\label{000}
The blue edges in $(X_1, B_1)$ forms a star.
\end{claim}

{\em Proof.} If there are $b_1, b_2\in B_1$ and $x_1, x_2\in X_1$ such that $b_1x_1, b_2x_2\in E^b(G)$. Note that there is a blue $b_1$ to $b_2$ path with at least three vertices in $(A_1, B_1)$. By Fact \ref{c7}, $x_1^+x_2^{++}, x_1^{--}x_2^{-}\in E^r(G)$ and $x_1^{--}b_2, x_2^{++}b_1\in E^r(G)$. Note $x_1^+, x_2^-\in Y$ and $x_1^{--}, x_2^{++}\in X_1\cup X_2$ and $x_1^{--}x_2^-, x_1^{--}b_2$ and $x_2^{++}x_1^+, x_2^{++}b_1$ are red edges and $(A_1, Y_2\cup B_2)$ and $(B_1, X_2\cup A_2)$ are red complete bipartite graphs. If $x_1^{--}, x_2^{++}\in X_2$, then by Lemma \ref{c3}, there is a red $C_{2m}$. If $x_1^{--}, x_2^{++}\in X_1$, then there are red paths with endpoints $x_1^+$ and $x_2^-$ in $(A_1, Y_2\cup B_2)$ of order from 3 to $min\{2|Y_2\cup B_2|-1, 2|A_1|+1\}$ and there are red paths with endpoints $b_1$ and $b_2$ in $(B_1, X_2\cup A_2)$ of order from 3 to $min\{2|X_2\cup A_2|+1, 2|B_1|-1\}$. Combining two of these paths with $x_1^{--}$ and $x_2^{++}$, there is a red $C_{2m}$. If $x_1^{--}\in X_1$ and $ x_2^{++}\in X_2$ (By symmetry, the case $x_1^{--}\in X_2$ and $ x_2^{++}\in X_1$ is similar), then there are red paths with endpoints $x_1^+$ and $x_2^-$ in $(A_1, Y_2\cup B_2)$ of order from 3 to $min\{2|Y_2\cup B_2|-1, 2|A_1|+1\}$ and there are red paths with endpoints $b_2$ and $x_2^{++}$ in $(B_1, X_2\cup A_2)$ of order from 2 to $min\{2|X_2\cup A_2|, 2|B_1|\}$. Combining two of these paths $P_1$, $P_2$ and $x_1^{--}$, we have a red $C_{2m}$: $x_1^{--}\rightarrow x_2^-$ (along $P_1$)$\rightarrow$  $x_1^+\rightarrow x_2^{++}$ (along $P_2$)$\rightarrow b_2\rightarrow x_1^{--}$. A contradiction. This complete the proof of Claim \ref{000}.\q

{\em Case 3.1.} The blue star in $(X_1, B_1)$ has center in $X_1$, i.e. $|X_1|=1$. Without loss of generality, assume that $X_1=\{x\}$.
\begin{claim}\label{00}
$(A_1, B_1)$ forms a blue complete bipartite graph.
\end{claim}
{\em Proof.}
Recall that $(A_1, B_2\cup Y_2)$ and $(B_1, A_2\cup X_2)$ are red complete bipartite graphs. If there is a red edge in $(A_1, B_1)$, by Lemma \ref{c3}, $(X_2\cup A_2, Y_2\cup B_2)$ forms a blue complete bipartite graph, otherwise there is a red $C_{2m}$. Since there is no blue $K_{n-1, n-1}$, $|X_2\cup A_2|\le n-2$, so $|A_1\cup X_1|=|A_1|+1\ge m+1$ if $n\ge m+1$, and $|A_1|+1\ge m+2$ if $n=m$. If $n\ge m+1$, since $(A_1, Y_2\cup B_2)$ forms a red complete bipartite, $G$ has a red $C_{2m}$. If $n=m$, and if $|Y_2|+|B_2|\ge m$, then $G$ contains a red $C_{2m}$ in $(A_1, Y_2\cup B_2)$. If $n=m$, and if $|Y_2|+|B_2|=m-1$, then $|B_1|= m+1$. Therefore $|A_1|, |B_1|\ge m+1$. By Lemma \ref{c3}, there is no disjoint red edges in $(A_1, B_1)$, then there is a blue $K_{m, m+1}$ in $(A_1, B_1)$. So $G$ contains a blue $C_{2n}$, a contradiction. This complete the proof of Claim \ref{00}.\q

If $n\ge m+1$, then $|A_1|\le m-1$ since otherwise there is a red $C_{2m}$ in $(A_1, Y_2)$. So $|A_1\cup\{x\}|\le m$. Then $|X_2\cup A_2|\ge n-1$. Since there is no red $C_{2m}$ in $(X_2\cup A_2, B_1)$, $|B_1|\le m-1$ and $|Y_2\cup B_2|\ge n$. We claim that $|A_1|<m-1$, otherwise $(A_1, Y_2\cup B_2)$ contains a red $K_{m-1, n}$. Since $G$ contains no red $C_{2m}$, each vertex in $X\setminus A_1$ has at most one red neighbor in $Y_2\cup B_2$. Since $G$ has no blue $C_{2n}$ and $|X\setminus A_1|, |Y_2\cup B_2|\ge n$, by Lemma \ref{c4}, there is a blue $K_{n-1, n-1}$ in $(X\setminus A_1, Y_2\cup B_2)$, a contradiction to Lemma \ref{c1}.
So we have shown that $|A_1|\le m-2$, consequently $|X_2\cup A_2|\ge n.$ Since there is no disjoint red edges in $(X_2\cup A_2, Y_2\cup B_2)$, there is a blue $K_{n-1, n-1}$ in $(X_2\cup A_2, Y_2\cup B_2)$, a contradiction to Lemma \ref{c1}.

If $m=n$, recall that $(A_1, B_1)$ forms blue complete bipartite graph. If $|B_1|\ge m$, then $|A_1|\le m-2$, otherwise $(A_1, B_1)$ contains a blue $K_{m-1, m-1}$, a contradiction to Lemma \ref{c1}. However if $|A_1|\le m-2$, then $|X_2\cup A_2|\ge m+1$, therefore $(X_2\cup A_2, B_1)$ contains a red $C_{2m}$, a contradiction. If $|B_1|\le m-1$, then $|Y_2\cup B_2|\ge m+1$, then $|A_1|\le m-1$ since otherwise the complete red bipartite graph $(A_1, Y_2\cup B_2)$ contains a red $C_{2m}$. So $|X_2\cup A_2|\ge m$. By Lemma \ref{c3}, there is no red disjoint edges in $(X_2\cup A_2, Y_2\cup B_2)$. So there is a blue $K_{m-1, m-1}$ in $(X_2\cup A_2, Y_2\cup B_2)$, a contradiction to Lemma \ref{c1}.

{\em Case 3.2.} The center of the blue star in $(X_1, B_1)$ is in $B_1$. Since we have shown for the case $|X_1|=1$, we may assume that $|X_1|\ge 2$ and there is $b\in B_1$ such that $(X_1, b)$ is a blue star.

\begin{claim}\label{0000}
$|A_1|+|B_1|=m$
\end{claim}
{\em Proof.}
Suppose on the contrary that $|A_1|+|B_1|\ge m+1$. We claim that $|A_1|\le m-1$. Recall that $(A_1, Y_2\cup B_2)$ forms a red complete bipartite graph and $|Y_2|=n-1\ge m-1$.   If $n=m$, by Lemma \ref{c1}, there is no red $K_{m-1, m-1}$ in $G$, so $|A_1|\le m-2$. If $n\ge m+1$, then $|Y_2|\ge m$. Since $(A_1, Y_2)$ is a red complete bipartite graph and has no red $C_{2m}$, $|A_1|\le m-1$. Note that $(B_1\setminus\{b\}, X\setminus A_1)$ and $(A_1, Y\setminus B_1)$ are red complete bipartite graphs. Note $|X\setminus A_1|\ge n$. By Lemma \ref{c3} and Lemma \ref{c1}, there is no disjoint red edges and no blue $K_{n-1, n-1}$ in $(X\setminus A_1, Y\setminus B_1)$. Then $|Y\setminus B_1|\le n-1$ and $|B_1|\ge m$ and $|B_1\setminus\{b\}|\ge m-1$ and $(B_1\setminus\{b\}, X\setminus A_1)$ contains a red $K_{m-1, n}$. Since $G$ has no red $C_{2m}$, each vertex form $Y_2\cup B_2\cup \{b\}$ has at most one red neighbor in $X\setminus A_1$. Since $|Y_2\cup B_2\cup \{b\}|\ge n$ and $|X\setminus A_1|\ge n$ and $G$ has no blue $C_{2n}$, by Lemma \ref{c4}, $(X\setminus A_1, Y_2\cup B_2\cup \{b\})$ contains a blue $K_{n-1, n-1}$, a contradiction to Lemma \ref{c1}. This completes the proof of Claim \ref{0000}.\q

If $n=m$, then $|A_1|+|B_1|\ge |A|=m+1$. By Claim \ref{0000}, we have $|A_1|+|B_1|=m$ and $n\ge m+1$. So $|A_1|=|B_2|$ and $|A_2|=|B_1|$. Since $|X\setminus A_1|\ge n$ and $|Y\setminus B_1|\ge n$ and there is no blue $C_{2n}$ and blue $K_{n-1, n-1}$ in $(X\setminus A_1, Y\setminus B_1)$, there are two disjoint red edges in $(X\setminus A_1, Y\setminus B_1)$. Note that $|X\setminus A_1|\ge n>|B_1\setminus\{b\}|$, and $|Y\setminus B_1|\ge n>|A_1|$, and $(B_1\setminus\{b\}, X\setminus A_1)$ and $(A_1, Y\setminus B_1)$ form red complete bipartite graphs. By Lemma \ref{c3}, there are all red cycles of length from $6$ to $2|A_1|+1+2|B_1\setminus\{b\}|+1=2m$, a contradiction. \q

\bigskip

\section{Proofs of Lemmas \ref{c1} and \ref{c2}}\label{sec3}

\subsection{Proof of Lemma \ref{c1}}

We first prove the following Lemma to be used in the proof of Lemma \ref{c1}.

\begin{lemma}\label{c6}
Let $G(A, B)$ and $H(C, D)$ be two bipartite graphs, where $A, B, C, D$ are pairwisely disjoint. Assume that $d_B(a)\ge |B|-1\ge 3$ and $d_D(c)\ge |D|-1\ge 3$ for all $a\in A$ and $c\in C$. If there are two disjoint edges between $A$ and $C$, then $G$ contains all even cycles of length $6$ to $min\{2|A|, 2|B|\}+min\{2|C|, 2|D|\}-2$.
\end{lemma}

In the proof of Lemma \ref{c6}, we will apply the following result by Jackson \cite{Jackson}.

\begin{theo}\label{jackson}(\cite{Jackson})
Let $G(X, Y)$ be a bipartite graph with $|X|\ge|Y|$. If $d(y)\ge max\{|Y|, \frac{|X|}{2}+1\}$ for every vertex $y\in Y$, then $G$ has a cycle containing all vertices in $Y$.
\end{theo}

{\em Proof of Lemma \ref{c6}.} Without loss of generality, we may assume that the disjoint edges are $a_1c_1$ and $a_2c_2$ with $a_1, a_2\in A$ and $c_1, c_2\in C$. Note that $d_B(a)\ge |B|-1\ge \frac{|B|}{2}+1$. By Theorem \ref{jackson}, $G(A\setminus\{a_2\}, B)$ contains all even cycle of length $4$ to $min\{2(|A|-1), 2(|B|-1)\}$ which contains $a_1$. Similarly, $H(C\setminus\{c_2\}, D)$ contains all even cycle of length $4$ to $min\{2(|C|-1), 2(|D|-1)\}$ which containing $c_1$.
Since $d_B(a_2)\ge |B|-1$, for any cycle in $G(A\setminus\{a_2\}, B)$ containing $a_1$, one of $a_1^-a_2$ and $a_1^+a_2$ must be contained in $E(G)$. Similarly, for any cycle in $H(C\setminus\{c_2\}, D)$ containing $c_1$, one of $c_1^-c_2$ and $c_1^+c_2$ must be contained in $E(H)$. Without loss of generality, we may assume that $a_1^+a_2$ and $c_1^+c_2$ are edges in $E(G)$ and $E(H)$. Combine this two cycles by using $a_1^+a_2, c_1^+c_2, a_1c_1, a_2c_2$ to replace $a_1a_1^+$ and $c_1c_1^+$, we get all even cycles of length $6$ to $min\{2(|A|-1), 2(|B|-1)\}+min\{2(|C|-1), 2(|D|-1)\}+2=min\{2|A|, 2|B|\}+min\{2|C|, 2|D|\}-2$.
\q

\noindent{\bf Lemma \ref{c1}} \
Let $N=m+n$ for $m=n\ge 5$ and $N=m+n-1$ for $n\ge m+1\ge 5$. If there is a red/blue-edge-coloring of $G(X, Y)=K_{N, N}$ such that there is neither a blue $C_{2n}$ nor a red $C_{2m}$. Then there is no blue $K_{n-1, n-1}$ in $G$.

\medskip

{\em Proof of Lemma \ref{c1}.} Suppose that $G$ has a blue $K_{n-1, n-1}$ with vertex set $X_1$ and $Y_1$ and $X_1\subsetneq X$ and $Y_1\subsetneq Y$ and $|X_1|=|Y_1|=n-1$. Let $A=X\setminus X_1$ and $B=Y\setminus Y_1$. Note that $|A|=|B|\ge m$. Let $A_1=\{x\in A| d_{Y_1}^b(x)\ge2\}$ and $X_1'=X_1\cup A_1$ and $B_1=\{y\in B| d_{X_1'}^b(y)\ge2\}$. Note that if both $|A_1|\ge 1$ and $|B_1|\ge 1$, then there is a blue $C_{2n}$ in $(X_1\cup A_1, Y_1\cup B_1)$. So we may always assume that $B_1=\emptyset$. Therefore $d_{X_1'}^b(y)\le 1$, i.e. $d_{X_1'}^r(y)\ge |X_1'|-1$ for all $y\in B$. If $|X_1'|-1\ge m$, then by Theorem \ref{jackson}, there is a red cycle $C_{2m}$. Therefore we only discuss the case $|X_1'|-1\le m-1$. Since $n\ge m$, $|A_1|=|X_1'|-(n-1)\le 1$, and if $|A_1|=1$, then $n=m$. So we have the following cases.

$|A_1|=0$ and $n=m+1$. Note that $d_{Y_1}^r(x)\ge n-2=m-1$ and $d_{X_1}^r(y)\ge n-2=m-1$ for $x\in A$ and $y\in B$. If $a_1w, a_2w\in E^r(G)$ for $a_1, a_2\in A$ and $w\in B$. By Theorem \ref{jackson}, there is a red cycle $C_{2m-2}$ in $(A\setminus a_2, Y_1)$ containing $a_1$. Note that either $a_2a_1^+$ or $a_2a_1^-$ ($a_1^+$ ($a_1^-$) is the successor (predecessor) in this red cycle) is red since $d_{Y_1}^r(a_2)\ge n-2=m-1$. Without loss of generality, assume $a_2a_1^+$ is red. We can replace $a_1a_1^+$ by $a_1w, wa_2, a_2a_1^+$ in the red cycle to get a red $C_{2m}$, a contradiction. So for each $y\in B$, $d_{A}^r(y)\le 1$ and $d_{A}^b(y)\ge |A|-1$. Similarly, for each $x\in A$, $d_{B}^r(x)\le 1$ and $d_{B}^b(x)\ge |B|-1$. By Lemma \ref{c6}, there is no disjoint red edges in $(A, B)$, recall that there is no red $K_{1, 2}$ in $(A, B)$, then $(A, B)$ contains a blue $K_{m, m}$ minus one edge. If there is a blue edge in $(X_1, B)$, then by Lemma \ref{c6}$, (A, Y_1)$ must form a red complete bipartite graph, and there is a red $C_{2m}$ in $(A, Y_1)$. If there is no blue edge in $(X_1, B)$, then $(X_1, B)$ forms a red complete bipartite graph, and there is a red $C_{2m}$ in $(X_1, B)$.

$|A_1|=0$ and $n=m$. In this case we have $|A|=|B|=m+1$. Note that $d_{Y_1}^r(x)\ge n-2=m-2$ and $d_{X_1}^r(y)\ge n-2=m-2$ for $x\in A$ and $y\in B$.
By Lemma \ref{c6}, there is no disjoint red edges in $(A, B)$, then there is a blue $K_{m, m}$ in $(A, B)$, a contradiction to that $G$ cannot have a blue $C_{2n}$.

$|A_1|=1$ and $n=m$. In this case we have $|A|=|B|=m+1$. Assume that $A'=A\setminus A_1$ and $A_1=\{a_1\}$. Then $d_{Y_1}^r(x)\ge n-2=m-2$ and $d_{X_1}^r(y)\ge n-2=m-2$ for $x\in A'$ and $y\in B$. By Lemma \ref{c6}, there is no disjoint red edges in $(A', B)$. Since $|A'|=m$ and $|B|=m+1$, $G$ has no blue $C_{2n}$, there is a blue $K_{m-1, m+1}=(A'', B)\subseteq (A', B)$. Since $G$ has no blue $C_{2n}$, each vertex in $X\setminus A''$ has at most one blue neighbor in $B$, i.e. each vertex in $X\setminus A''$ has at least $m$ red neighbors in $B$. By Theorem \ref{jackson}, there is a red $C_{2m}$ in $(X\setminus A'', B)$, a contradiction.\q

\subsection{Proof of Lemma \ref{c2}}

To prove Lemma \ref{c2}, we first prove Lemmas \ref{c4} and \ref{c9}. The following  results by Hu \cite{Huz}, and Li-Ning \cite{Lining} will be applied in the proofs of Lemmas \ref{c4} and \ref{c9}.

\begin{theo}\label{huz}(\cite{Huz})
Let $G(X, Y)$ be a hamiltonian bipartite graph on $2n$ vertices, where $n>3$. If $e(G)>\frac{n(n-1)}{2}+2$, then $G$ contains cycles of every possible even length.
\end{theo}

\begin{theo}\label{Lining}(\cite{Lining})
Let $G(X, Y)$ be a bipartite graph and $C$ be a longest cycle of $G$. Suppose that $|X|=n$, $|Y|=m$ and $|C|=2c$, where $n\ge m\ge c$. If $m\le 2c$, then $e(G-C)+e(G-C, C)\le c(n-1-c)+m$. If $m\ge 2c$, then $e(G-C)+e(G-C, C)\le c(m+n+1-3c)$.
\end{theo}

\begin{lemma}\label{c4}
Let $G(X, Y)$ be a bipartite graph with $|X|=|Y|\ge n\ge 4$. If $d(y)\ge |X|-1$ for every vertex $y\in Y$ and $G$ has no $C_{2n}$, then $K_{n-1, n-1}\subseteq G$.
\end{lemma}
{\em Proof of Lemma \ref{c4}.} If $|X|=|Y|\ge n+1$, then $d(y)\ge |X|-1\ge n$ for $y\in Y$. Take $Y'\subseteq Y$ with $|Y'|=n$ and apply Theorem \ref{jackson} to $(X, Y')$, then $C_{2n}\subseteq (X, Y')$, a contradiction. If $|X|=|Y|=n$, then $d(y)\ge n-1$ for $y\in Y$ and $e(G)\ge (n-1)n$. By Theorem \ref{jackson}, then $C=C_{2(n-1)}\subseteq G$. Since $G$ contains no $C_{2n}$, if $G$ has no $K_{n-1, n-1}$, by Theorem \ref{Lining}, we have $e(G)\le (n-1)(n-1)-1+n=(n-1)n$. So the $C_{2n-2}$ forms a $K_{n-1, n-1}$ minus one edges and there are $n$ edges in $E(G-C)\cup E(G-C, C)$. Let $x_0=X\setminus C$ and $y_0=Y\setminus C$. Note that $d(y_0)\ge n-1$. If $x_0y_0\in E(G)$, then $d_C(x_0)=0$, otherwise there is a $C_{2n}$ in $G$. Recall that there are $n$ edges in $E(G-C)\cup E(G-C, C)$, i.e. there are $n$ edges incident to $x_0$ or $y_0$, then $d_C(y_0)=n-1$. So $(X\cap C, Y\cap C\cup\{y_0\})$ forms a $K_{n-1, n}$ minus one edge. So there is a $K_{n-1, n-1}$ in $(X\cap C, Y\cap C \cup\{y_0\})$. If $x_0y_0\notin E(G)$, since $d(y_o)\ge n-1$, we have $d(y_0)=d_C(y_0)=n-1$. Therefore $(X\cap C, Y\cap C\cup\{y_0\})$ forms a $K_{n-1, n}$ minus one edge, there is a $K_{n-1, n-1}$ in $(X\cap C, Y\cap C \cup\{y_0\})$. \q

\begin{lemma}\label{c9}
Let $G(X, Y)$ be a bipartite graph with $|X|=|Y|=M$. Let $C$ be a longest cycle in $G$ with $|V(C)|=2c$. Assume that $G$ is $C_{2n}$-free. If $c\ge \frac{M}{2}$ and $c>n$, then $e(G)\le \frac{M^2-M+6}{2}$. If $c\le n-1 \le \frac{M}{2}$, then $e(G)\le (n-1)(2M-2n+3)$, and if equality holds, then $G[V(C)]$ must be a $K_{c, c}$.
\end{lemma}
{\em Proof of Lemma \ref{c9}.} If $c\ge \frac{M}{2}$ and $c>n$, by Theorem \ref{huz} and Theorem \ref{Lining}, we have
\begin{eqnarray}
e(G)&=&e(G[V(C)])+e(G-C, C)+e(G-C)\nonumber \\
&\le&\frac{c(c-1)}{2}+2+c(M-1-c)+M \nonumber \\
&=&-\frac{c^2}{2}+c(M-\frac{3}{2})+M+2 \nonumber\\
&\le&\frac{M^2-M+6}{2}.\nonumber
\end{eqnarray}
If $c\le n-1 \le \frac{M}{2}$, by Theorem \ref{Lining}, we have
\begin{eqnarray}
e(G)&=&e(G[V(C)])+e(G-C, C)+e(G-C)\nonumber \\
&\le&c^2+c(2M+1-3c) \nonumber \\
&=&-2c^2+c(2M+1) \nonumber\\
&\le&(n-1)(2M-2n+3).\nonumber
\end{eqnarray}
\q

\noindent{\bf Lemma \ref{c2}} \
Let $N=m+n$ for $m=n$ and $N=m+n-1$ for $n\ge m+1$. Assume that there is a red/blue-edge-coloring of $G(X, Y)=K_{N, N}$ such that there is neither a blue $C_{2n}$ nor a red $C_{2m}$. Assume that $C$ is a longest blue cycle with $|V(C)|=2c$ and $C'$ is a longest red cycle with $|V(C')|=2t$. If $n\ge m+1$, then $c\le n-1$. If $n=m$, then $c\le m-1$ and $t\le m-1$.

{\em Proof of Lemma \ref{c2}.} {\em Case 1.} $m=n$, $N=2m$.

Without loss of generality, assume that $c\ge m+1$. Then by Lemma \ref{c9}, $e^b(G)\le \frac{N^2-N+6}{2}$. If $t\ge m+1$, then by Lemma \ref{c9}, $e^r(G)\le \frac{N^2-N+6}{2}$, a contradiction to that $e^b(G)+e^r(G)=N^2$. If $t\le m-1$, then by Lemma \ref{c9}, $e^r(G)\le (m-1)(2N-2m+3)=(m-1)(2m+3)$. Note that $\frac{N^2-N+6}{2}+(m-1)(2m+3)=N^2$, then $e^r(G)=(m-1)(2m+3)$. So $t=m-1$ and $G[V(C')]$ forms a red $K_{m-1, m-1}$ (Lemma \ref{c9}), a contradiction to Lemma \ref{c1}.

{\em Case 2.} $n\ge m+1$, and $N=n+m-1\ge 2m$.

Suppose that $c\ge n+1$, then by Lemma \ref{c9}, $e^b(G)\le \frac{N^2-N+6}{2}$. Therefore $e^r(G)\ge\frac{N^2+N-6}{2}$.

If $t\ge \frac{N}{2}$, then by Lemma \ref{c9}, $e^r(G)\le \frac{N^2-N+6}{2}$, a contradiction.

If $m+1\le t\le \frac{N}{2}$, then by Theorem \ref{huz} and Theorem \ref{Lining}, we have
\begin{eqnarray}
e^r(G)&\le&\frac{t(t-1)}{2}+2+t(2N+1-3t) \nonumber \\
&=&-\frac{5t^2}{2}+t(2N+\frac{1}{2})+2 \nonumber\\
&\le&\frac{4N^2+2N+20.25}{10}.\nonumber
\end{eqnarray}
Then $\frac{4N^2+2N+20.25}{10}\ge\frac{N^2+N-6}{2}$, i.e. $N\le \frac{\sqrt{210}-3}{2}< 6$, a contradiction.

If $t\le m-1$, then by Lemma \ref{c9}, $e^r(G)\le(m-1)(2N-2m+3)$. Therefore $(m-1)(2N-2m+3)\ge\frac{N^2+N-6}{2}$, i.e. $\frac{N^2}{2}-(2m-\frac{5}{2})N+2m^2-5m\le 0$. We have $2m-5\le N\le 2m$. Since $N\ge 2m$, we have $N=2m$ and $n=m+1$. Since the equality holds, $e^r(G)=(m-1)(2m+3)$ and $C'$ forms a red $K_{m-1, m-1}$. Without loss of generality, assume that $(X_1, Y_1)$ forms a red $K_{m-1, m-1}$ with $X_1\subset X$ and $Y_1\subseteq Y$ and $|X_1|=|Y_1|=m-1$. Let $A=X\setminus X_1$ and $B=Y\setminus Y_1$. Note that $|A|=|B|= m+1$. Let $A_1=\{x\in A| d_{Y_1}^r(x)\ge2\}$ and $X_1'=X_1\cup A_1$ and $B_1=\{y\in B| d_{X_1'}^r(y)\ge2\}$. Since $G$ has no red $C_{2m}$, we may assume that $B_1=\emptyset.$ Note that each vertex in $B$ has at most one red neighbor in $X_1'$. If $|A_1|\ge 2$, then $|X_1'|\ge|B|=m+1=n$. Since $G$ has no blue $C_{2n}$, by Lemma \ref{c4}, we get a blue $K_{n-1, n-1}$ in $(X_1', B)$, a contradiction to Lemma \ref{c1}. Therefore $|X_1'|\le m$. By Lemma \ref{c6}, there is no disjoint blue edges in $(X\setminus X_1', B)$. 
If $|A_1|=0$, then $(X\setminus X_1', B)$ contains a red $K_{m, m+1}$, a contradiction. If $|A_1|=1$, since $G$ has no red $C_{2m}$, $(X\setminus X_1', B)$ contains a red $K_{m-1, m+1}=(A', B)$ with $A'\subset X\setminus X_1'$. So each vertex in $X\setminus A'$ has at most one red neighbor in $B$. Note that $|X\setminus A'|=|B|=m+1=n$. Since $G$ has no blue $C_{2n}$, by Lemma \ref{c4}, there is a blue $K_{n-1, n-1}$ in $(X\setminus A', B)$, a contradiction to Lemma \ref{c1}.
\q

\bigskip

{\bf Remarks.} For more than two colors, we have very few results on exact values of Ramsey numbers. A nice result due to Jenssen and Skokan \cite{JS} gives $r^k(C_{2n+1})$ (a conjecture of Bondy and Erd\H os)  for large enough $n$. Most current known results are in diagonal cases (the same graph for all colors). For non-diagonal cases, we even do not know $r(P_l, P_m, P_n)$ and $br(P_l, P_m, P_n)$.

\end{document}